\documentclass{elsart}
\usepackage{amssymb,amsfonts}
\newcommand{\LHS}{\mbox{LHS}}
\newcommand{\RHS}{\mbox{RHS}}
\newcommand{\trace}{\mathop{\rm Tr}\nolimits}

\newcommand{\re}{\mathop{\rm Re}\nolimits}

\newcommand{\real}{\re}

\newcommand{\twomat}[4]{\left(\begin{array}{cc}#1&#2\\#3&#4\end{array}\right)}

\newcommand{\schatten}[2]{\left|\left|\,{#2}\,\right|\right|_{#1}}

\newcommand{\cA}{{\mathcal A}}
\newcommand{\cB}{{\mathcal B}}
\newcommand{\cC}{{\mathcal C}}
\newcommand{\cH}{{\mathcal H}}

\newcommand{\C}{{\mathbb{C}}}
\newcommand{\R}{{\mathbb{R}}}

\DeclareRobustCommand\openone{\leavevmode\hbox{\small1\normalsize\kern-.33em1}}
\newcommand{\id}{\mathrm{\openone}}
\newcommand{\identity}{\id}
\newcommand{\ddt}{\left.\frac{\partial}{\partial t}\right|_{t=0}}

\newcommand{\be}{\begin{equation}}
\newcommand{\ee}{\end{equation}}
\newcommand{\bea}{\begin{eqnarray}}
\newcommand{\eea}{\end{eqnarray}}
\newcommand{\beas}{\begin{eqnarray*}}
\newcommand{\eeas}{\end{eqnarray*}}

\newtheorem{lemma}{Lemma}
\newtheorem{corollary}{Corollary}

\newtheorem{proposition}{Proposition}
\newtheorem{conjecture}{Conjecture}
\newcount\minute
\newcount\hour
\def\currenttime{%
    \minute\time
    \hour\minute
    \divide\hour60
    \the\hour:\multiply\hour60\advance\minute-\hour\the\minute}
\begin{document}
\begin{frontmatter}
\title{On a Norm Compression Inequality for $2\times N$ Partitioned Block Matrices}
\author{Koenraad M.R. Audenaert}
\address{
Institute for Mathematical Sciences, Imperial College London,\\
53 Prince's Gate, London SW7 2PG, United Kingdom}
\ead{k.audenaert@imperial.ac.uk}
\date{\today, \currenttime}
\begin{abstract}
We conjecture the following so-called norm compression inequality for
$2\times N$ partitioned block matrices and the Schatten $p$-norms:
for $p\ge 2$,
$$
\schatten{p}{
\left(
\begin{array}{cccc}
A_1 & A_2 & \cdots & A_N \\
B_1 & B_2 & \cdots & B_N
\end{array}
\right)}
 \le \schatten{p}{
\left(
\begin{array}{cccc}
||A_1||_p & ||A_2||_p & \cdots & ||A_N||_p \\
||B_1||_p & ||B_2||_p & \cdots & ||B_N||_p
\end{array}
\right)}
$$
while for $1\le p\le 2$ the ordering of the inequality is reversed.
This inequality includes Hanner's inequality for matrices as a special case.
We prove several special cases of this inequality and give examples
for $3\times 3$ and larger partitionings where it does not hold.
\end{abstract}
\end{frontmatter}
\section{Introduction\label{sec1}}
The \textit{norm compression} of a block-partitioned matrix $T=[T_{(ij)}]$ w.r.t.\ a given matrix norm $||.||$ is a matrix obtained
from $T$ by replacing each of its blocks by their norm: $[||T_{(ij)}||]$.
One can raise the question as to how the norm of $T$ relates to the norm of its norm compression.
In many cases one can find simple upper and lower bounds to one in terms of the other, giving
rise to \textit{norm compression inequalities} (NCI's); see \cite{ka} and references
therein for an overview of several of these inequalities.

The Schatten $p$-norms, for $1\le p<\infty$, are unitarily invariant (UI) norms
generalising the $l_p$ norms to the non-commutative setting.
For a general matrix or operator $A$, they are defined as
$$
||A||_p = (\trace(|A|^p))^{1/p},
$$
which reduces for positive semi-definite (PSD) matrices $A$ to
$$
||A||_p = (\trace(A^p))^{1/p}.
$$
If one keeps this definition as is and considers the case $0\le p\le 1$, one obtains quasi-norms.
The $p$-norms are convex for $p\ge1$; this is just the triangle inequality.
For $0\le p\le 1$, the $p$-quasi-norms are concave on the positive semi-definite (PSD) matrices;
this is no longer true when dealing with non-PSD matrices, since the absolute value operation then comes into play,
which is not concave.
To prove this convexity/concavity statement, one can for example apply Theorem 2 in \cite{carlenlieb}
to a matrix $A$ on $\C^2\otimes\cH$ of the form $A=A_1\oplus A_2$, with $A_{1,2}\ge0$.

In this paper, I propose a norm compression inequality, based on the Schatten norms,
for $2\times N$ block-partitioned matrices:
$$
T = \left(
\begin{array}{cccc}
A_1 & A_2 & \cdots & A_N \\
B_1 & B_2 & \cdots & B_N
\end{array}
\right).
$$
I denote by $\cC_p(T)$ its Schatten $p$-norm compression,
$$
\cC_p(T) =
\left(
\begin{array}{cccc}
||A_1||_p & ||A_2||_p & \cdots & ||A_N||_p \\
||B_1||_p & ||B_2||_p & \cdots & ||B_N||_p
\end{array}
\right).
$$
I conjecture the following:
\begin{conjecture}\label{conj1}
Let $T$ be a general matrix partitioned in $2\times N$ blocks,
and let $\cC_p(T)$ be its norm compression using the Schatten $p$-norm,
then the following norm compression inequalities hold:
\bea
||T||_p &\ge& \schatten{p}{\cC_p(T)},\quad 1\le p\le2 \label{eq:thedual}\\
||T||_p &\le& \schatten{p}{\cC_p(T)},\quad p\ge2. \label{eq:the}
\eea
\end{conjecture}
The work presented in this paper grew out of on-going work to prove Hanner's inequality \cite{hanner} for matrices.
Hanner's inequality proper \cite{hanner} for the $L_p$ function spaces is, for $f,g$ functions in $L_p$,
\be
||f+g||_p^p + ||f-g||_p^p \ge (||f||_p+||g||_p)^p + \big|\, ||f||_p-||f||_p \big|^p,
\ee
for $1\le p\le 2$, and the reversed inequality for $2\le p$.
It is widely believed that these inequalities are also true for the Schatten trace ideals $\cC_p$:
for general matrices $A,B$ in $\cC_p$, and $1\le p\le 2$,
\be\label{hanner}
||A+B||_p^p + ||A-B||_p^p \ge (||A||_p+||B||_p)^p + \big|\, ||A||_p-||B||_p \big|^p,
\ee
while for $2\le p$, the inequality is reversed.
This generalisation to matrices has been proven in a number of instances:
\begin{enumerate}
\item For $1\le p$, when $A+B$ and $A-B$ are PSD.
\item For $1\le p\le4/3$, $p=2$, and $4\le p\le\infty$, when $A$ and $B$ are general matrices.
\end{enumerate}
Proofs are due to Ball, Carlen and Lieb \cite{lieb} and Tomczak-Jaegermann \cite{TJ}.

That Conjecture \ref{conj1} implies Hanner's inequality for matrices can be seen easily
by putting $N=2$, $A_1=B_2=A$ and $A_2=B_1=B$.
Unitarily conjugating
$T$ with the matrix $\twomat{\id}{\id}{\id}{-\id}$, and its norm compression with
$\twomat{1}{1}{1}{-1}$ directly yields (\ref{hanner}).

Two other special cases where Conjecture \ref{conj1} is proven are for $2\times 2$-partitioned block matrices $T$ ($N=2$).
The case when $T$ is PSD has been proven by King \cite{king},
and the case when the blocks of $T$ are all diagonal matrices has been proven by King and Nathanson \cite{king_nath}.

Note that, if true, (\ref{eq:the}) is strongly sharp,
which means that equality can be obtained for any given imposed values $a_k\ge0$ and
$b_k\ge0$ on $||A_k||_p$ and $||B_k||_p$, respectively.
To see this, one just takes $A_k=[a_k]\oplus 0$ and $B_k=[b_k]\oplus 0$.

In the following, I present equivalent forms of the Conjecture and a number of proofs in some further special cases.
Section 2 deals with a duality argument to show that the $2\le p$ case can be derived from the $1\le p\le 2$ case,
while the converse is, surprisingly, not clear at all.
Equivalent forms of the Conjecture are given in Section 3.
The central part of the paper is Section 4, in which I give proofs of the Conjecture for 4 special cases:
\begin{itemize}
\item The norm compression of $T$ has rank 1.
\item All blocks in $T$ have rank 1.
\item All blocks $A_k$ in the first row are proportional, and so are all blocks $B_k$ in the second row.
\item Just as in the case of Hanner's inequality, $p\ge4$.
\end{itemize}
For the first three special cases, it turns out that (\ref{eq:thedual}) also holds for $0\le p\le 1$.
This, however, is not generally true.

Given the fact that Conjecture \ref{conj1} is a generalisation of Hanner's inequality to $2\times N$ block
matrices, one may wonder whether the Conjecture could even be true for arbitrary $M\times N$ partitionings.
It turns out that this is not the case.
For $4\times 4$ block matrices there are counterexamples when the blocks $A_{ij}$ are real scalars,
in which case the norm-compression is just the entry-wise absolute value $[|A_{ij}|]$.
The following example already appeared in \cite{ka}: for
$$
A=\left(\begin{array}{rrrr}
     3  &   0 &   -2  &  -2 \\
     0  &   3 &    2  &  -1 \\
    -2  &   2 &    4  &   0 \\
    -2  &  -1 &    0  &   3
\end{array}
\right)\ge 0,
$$
one finds $\schatten{1.5}{A}=9.49929...$ and $\schatten{1.5}{[|A_{ij}|]}=9.63184...$.

In the final Section of this paper I consider the special case of diagonal blocks, and find
counterexamples for $3\times 3$ block matrices, and counterexamples to the extension
of (\ref{eq:thedual}) to $0\le p\le 1$.
\section{Duality}
A standard duality argument allows to reduce the amount of effort in proving Conjecture \ref{conj1} by a factor of two.
I will show that validity of (\ref{eq:thedual}) implies validity of (\ref{eq:the}).
Although seemingly obvious, proving the converse implication turns out to be problematic.
\subsection{(\ref{eq:thedual}) implies (\ref{eq:the})}
Let $1\le p\le 2$ and let $p'$ be its conjugate exponent $1/p+1/p'=1$; thus $p'\ge2$.

Using the dual representation of a norm, we can express the left-hand side (LHS) of (\ref{eq:the})
in a variational manner as
$$
\schatten{p'}{\twomat{A_1}{\ldots}{B_1}{\ldots}} \\
= \max_{S_k, T_k} \left\{ \real\trace
\twomat{A_1}{\ldots}{B_1}{\ldots}
\twomat{S_1}{\ldots}{T_1}{\ldots}^* :
\schatten{p}{
\twomat{S_1}{\ldots}{T_1}{\ldots}
}\le 1
\right\},
$$
where the maximisation is over all real scalars $s_1,\ldots,s_N$ and $t_1,\ldots,t_N$.
Thus,
$$
||T||_{p'} = \trace[T\hat{Y}^*],
$$
where $\hat{Y}=\twomat{S_1}{\ldots}{T_1}{\ldots}$ is the optimal $Y$ in the above maximisation, with $\schatten{p}{\hat{Y}}=1$.

By the assumed validity of (\ref{eq:thedual}), we have
$$
\schatten{p}{\twomat{||S_1||_p}{\ldots}{||T_1||_p}{\ldots}}\le 1.
$$
The RHS of (\ref{eq:the}) can be expressed as
\beas
\RHS(\ref{eq:the}) &=&
\schatten{p'}{\twomat{||A_k||_{p'}}{\ldots}{||B_k||_{p'}}{\ldots}} \\
&=& \max_{s_k, t_k} \left\{
\sum_{k=1}^N s_k ||A_k||_{p'} + t_k ||B_k||_{p'} :
\schatten{p}{
\twomat{s_1}{\ldots}{t_1}{\ldots}
}\le 1
\right\},
\eeas
where the maximisation is over all real scalars $s_1,\ldots,s_N$ and $t_1,\ldots,t_N$.
Likewise, we have
\beas
||A_k||_{p'} &=& \max_{G_k} \{\real\trace[A_k G_k^*] : ||G_k||_p=1 \} \\
||B_k||_{p'} &=& \max_{H_k} \{\real\trace[B_k H_k^*] : ||H_k||_p=1 \}.
\eeas
Thus we obtain
\beas
\RHS(\ref{eq:the})
&=& \max_{s_k, t_k, G_k, H_k}
\real\sum_{k=1}^N s_k \trace[A_k G_k^*] + t_k \trace[B_k H_k^*],
\eeas
under the aforementioned constraints.

Now note that if we put
$s_k = ||\hat{S}_k||_p$, $t_k=||\hat{T}_k||_p$, $G_k=\hat{S}_k/s_k$ and $H_k=\hat{T}_k/t_k$,
we obtain $s_k,t_k,G_k,H_k$ satisfying the constraints of the above maximisation.
Furthermore, with these values we get the identity
\beas
||T||_{p'} = \trace[T\hat{Y}^*] &=& \sum_k \trace[A_k \hat{S}_k^*+B_k \hat{T}_k^*]\\
&=&\sum_k \trace[A_k s_k G_k^* + B_k t_k H_k^*]\\
&=&\sum_k s_k\trace[A_k G_k^*] + t_k\trace[B_k H_k^*],
\eeas
which is real for this particular choice.
We conclude that this expression is upper bounded by $\RHS(\ref{eq:the})$, which is the
statement of (\ref{eq:the}).
\qed
\subsection{Does (\ref{eq:the}) imply (\ref{eq:thedual})?}
One might think that it would be equally easy to prove the converse duality,
i.e. that validity of (\ref{eq:the}) implies validity of (\ref{eq:thedual}).
This turns out to be not so easy, because of an unexpected extra difficulty,
and we currently have no proof of any converse duality.

Let us again take $p$ between 1 and 2.
Using the dual representation of a norm, we again express the LHS of (\ref{eq:thedual})
in a variational manner as
\beas
\schatten{p}{T} &=&
\schatten{p}{\twomat{A_1}{\ldots}{B_1}{\ldots}} \\
&=& \max_{S_k, T_k} \left\{ \real\trace
\twomat{A_1}{\ldots}{B_1}{\ldots}
\twomat{S_1}{\ldots}{T_1}{\ldots}^* :
\schatten{p'}{
\twomat{S_1}{\ldots}{T_1}{\ldots}
}\le 1
\right\} \\
&=& \max_{S_k, T_k} \left\{
\real\sum_{k=1}^N \trace[A_k S_k^* + B_k T_k^*] :
\schatten{p'}{
\twomat{S_1}{\ldots}{T_1}{\ldots}
}\le 1
\right\},
\eeas
where the maximisation is over all blocks $S_1,\ldots,S_N$ and $T_1,\ldots,T_N$ of the same size as
$A_k$ and $B_k$.

Now consider RHS(\ref{eq:thedual}), which can similarly be expressed as
$$
\schatten{p}{\twomat{||A_k||_p}{\ldots}{||B_k||_p}{\ldots}}
= \max_{s_k, t_k} \left\{
\sum_{k=1}^N s_k ||A_k||_p + t_k ||B_k||_p :
\schatten{p'}{
\twomat{s_1}{\ldots}{t_1}{\ldots}
}= 1
\right\},
$$
where the maximisation is over all real scalars $s_1,\ldots,s_N$ and $t_1,\ldots,t_N$.
Let us denote the optimal values by $\hat{s}_k,\hat{t}_k$ and introduce the matrix
$$
\hat{X}=\twomat{\hat{s}_1}{\ldots}{\hat{t}_1}{\ldots}.
$$
Its entries by no means need to be positive, which will turn out to be a significant complication.

Likewise, we have
\beas
||A_k||_p &=& \max_{G_k} \{\real\trace[A_k G_k^*] : ||G_k||_{p'}=1 \} \\
||B_k||_p &=& \max_{H_k} \{\real\trace[B_k H_k^*] : ||H_k||_{p'}=1 \}.
\eeas
We denote the optimal $G_k$ and $H_k$ by $\hat{G}_k,\hat{H}_k$.

Inserting these expressions in the previous one gives
$$
\RHS(\ref{eq:thedual}) =
\sum_k \trace[A_k \, \hat{s}_k \hat{G}_k^* + B_k \,\hat{t}_k \hat{H}_k^*].
$$

Comparing this expression for RHS(\ref{eq:thedual}) to the one for LHS(\ref{eq:thedual}), we see
that $\LHS\ge \RHS$ follows if on putting $S_k=\hat{s}_k \hat{G}_k$ and $T_k=\hat{t}_k \hat{H}_k$ we obtain a matrix
$$
Y:=\twomat{\hat{s}_1 \hat{G}_1}{\ldots}{\hat{t}_1 \hat{H}_1}{\ldots}
$$
that satisfies the constraint
$\schatten{p'}{Y}\le 1$.

If the optimal $\hat{s}_k,\hat{t}_k$ happen to be non-negative, this is indeed the case.
Combining the constraints on $s_k, t_k, G_k, H_k$ with the assumption
that Conjecture \ref{conj1} holds for $p'$, we find
\beas
\schatten{p'}{Y}&\le&
\schatten{p'}{
\twomat{|\hat{s}_1|\,\, ||\hat{G}_1||_{p'}}{\ldots}{|\hat{t}_1| \,\,||\hat{H}_1||_{p'}}{\ldots}
}
=\schatten{p'}{
\twomat{\hat{s}_1}{\ldots}{\hat{t}_1}{\ldots}
}= 1.
\eeas

Alas, this no longer works if some $\hat{s}_k,\hat{t}_k$ are negative, because then it can happen
that $\schatten{p'}{Y}$ is greater than 1
(the author is grateful to Chris King for pointing this out).

Take, for example, $p=1.5$ ($p'=3$) and the $2\times 2$ matrix $T$ with scalar blocks
$T=\twomat{1}{2}{-0.1}{1}$.
Its norm compression is $A:=\cC_p(T)=\twomat{1}{2}{0.1}{1}$, and the optimal $X$ is
$\hat{X}=\twomat{0.488}{0.774}{-0.153}{0.488}$ (which has 3-norm equal to 1),
while $\hat{G}_1=\hat{G}_2=\hat{H}_2=1$ and $\hat{H}_1=-1$.
This leads to
$$
Y=\twomat{\hat{s}_1 \hat{G}_1}{\ldots}{\hat{t}_1 \hat{H}_1}{\ldots} = \twomat{0.488}{0.774}{+0.153}{0.488},
$$
which has 3-norm 1.0426, i.e.\ larger than 1.

Nevertheless, $\RHS=\schatten{p}{A}=\trace[AX^*]=\trace[TY^*]<\trace[TX^*]$, because of the sign patterns,
and $\trace[TX^*]\le\schatten{p}{T}=\LHS$, because $X$ has 3-norm equal to 1. In this example,
the conjecture itself is not violated,
but the implication $(\ref{eq:the})\Longrightarrow(\ref{eq:thedual})$ does not go through.

\section{Equivalent and Related Forms}
Let us consider first the case $p\ge2$.
We put $p=2q$, with $q\ge 1$, and perform the conversion
$||T||_{2q} = ||TT^*||_q^{1/2}$.
Inequality (\ref{eq:the}) then becomes, after squaring both sides: for $q\ge 1$,
\be\label{eq:the2}
\schatten{q}{
\twomat{\sum_k A_k A_k^*}{\quad\sum_k A_k B_k^*}{\sum_k B_k A_k^*}{\quad\sum_k B_k B_k^*}
} \le \schatten{q}{
\twomat{\sum_k a_k^2}{\quad\sum_k a_k b_k}{\sum_k a_k b_k}{\quad\sum_k b_k^2}
},
\ee
with
$a_k = ||A_k||_{2q} = ||A_k A_k^*||_q^{1/2}$ and
$b_k = ||B_k||_{2q} = ||B_k B_k^*||_q^{1/2}$.
Note that both matrices appearing in (\ref{eq:the2}) are PSD.

In the following I will use the subscript $(jk)$ to indicate the $j,k$-th \textit{block},
as opposed to the subscript $jk$ without the brackets, which I use to denote the $j,k$-th \textit{entry}.

An equivalent formulation of (\ref{eq:the2}) is the following
\begin{conjecture}\label{conj3}
Let $Q^{(i)}=[Q^{(i)}_{(jk)}]_{j,k=1}^2$ be arbitrary PSD $2\times 2$ block matrices, and $q^{(i)}_k$ arbitrary non-negative numbers.
The maximum attainable value of $||\sum_i Q^{(i)}||_q$ under the constraints $||Q^{(i)}_{(kk)}||_q = q^{(i)}_k$
is given by
\bea
\lefteqn{\max_{Q^{(i)}\ge0} \{||\sum_i Q^{(i)}||_q: ||Q^{(i)}_{(kk)}||_q = q^{(i)}_k\}} \nonumber \\
&=&
\schatten{q}{
\twomat{\sum_i q^{(i)}_1}{\quad\sum_i \sqrt{q^{(i)}_1 q^{(i)}_2}}{\sum_i \sqrt{q^{(i)}_1 q^{(i)}_2}}{\quad\sum_i q^{(i)}_2}
}.\label{eq:QRbound2}
\eea
\end{conjecture}
\textit{Proof of equivalence of (\ref{eq:the2}) and (\ref{eq:QRbound2}).}
It is immediately clear that (\ref{eq:QRbound2}) implies (\ref{eq:the2}), as can be seen by setting
$Q^{(i)}:=\twomat{A_i A_i^*}{A_i B_i^*}{B_i A_i^*}{B_i B_i^*}$,
$q^{(i)}_1=a_i^2$ and $q^{(i)}_2=b_i^2$.

To prove the converse,
I first show that maximality of $||\sum_i Q^{(i)}||_q$ under the above-mentioned constraints
implies that all $Q^{(i)}$ are of the form
\be
Q^{(i)}=\twomat{A_i A_i^*}{A_i B_i^*}{B_i A_i^*}{B_i B_i^*}. \label{eq:Qreq2}
\ee

W.l.o.g.\ we can assume all blocks having the same size (if not, the smaller
blocks can be padded with zeroes, an operation which does not change the Schatten norms).
We can maximise $||\sum_i Q^{(i)}||_q$ in two steps.
We first maximise it keeping the diagonal blocks $Q^{(i)}_{(kk)}$ fixed,
and then maximise over the diagonal blocks while keeping their norms fixed.

In the first step, positivity of $Q^{(i)}$ requires
$Q^{(i)}_{(12)} = (Q^{(i)}_{(11)})^{1/2} K^{(i)} (Q^{(i)}_{(22)})^{1/2}$, where $K^{(i)}$ is a contraction.
Maximising $||\sum_i Q^{(i)}||_q$ over all allowed $Q^{(i)}$
then amounts to maximising it over all contractions $K^{(i)}$.
The set of square contractions is convex and its extremal points
are the unitaries.
Since $||.||_q$ is a convex function of its argument, $||\sum_i Q^{(i)}||_q$ will achieve its maximum
in an extremal point, i.e.\ in a point where
$Q^{(i)}_{(12)} = (Q^{(i)}_{(11)})^{1/2} U^{(i)} (Q^{(i)}_{(22)})^{1/2}$, for some unitary $U^{(i)}$.
Furthermore, this remains so when performing the final optimisation over the diagonal blocks.
This means that the optimal $Q^{(i)}$ are indeed of the required form (\ref{eq:Qreq2}),
with $A_i = (Q^{(i)}_{(11)})^{1/2}$, and $B_i = U^{(i)}(Q^{(i)}_{(22)})^{1/2}$.

Applying (\ref{eq:the2}) then immediately yields (\ref{eq:QRbound2}),
since with these assignments we have
$q^{(i)}_1=a_i^2$, and $q^{(i)}_2=b_i^2$.
\qed


\bigskip

A (potentially) more general statement than (\ref{eq:the2}) is the following
\begin{conjecture}\label{conj2}
Let $Q=[Q_{(jk)}]$ and $R=[R_{(jk)}]$ be arbitrary PSD block matrices, and $q_k, r_k$ arbitrary non-negative numbers.
The maximum attainable value of $||Q+R||_q$ under the constraints $||Q_{(kk)}||_q = q_k$ and $||R_{(kk)}||_q = r_k$
is given by
\bea
\lefteqn{\max_{Q,R\ge0} \{||Q+R||_q: ||Q_{(kk)}||_q = q_k, ||R_{(kk)}||_q = r_k\}} \nonumber \\
&=&
\schatten{q}{
\twomat{\sum_k q_k}{\quad\sum_k \sqrt{q_k r_k}}{\sum_k \sqrt{q_k r_k}}{\quad\sum_k r_k}
}.\label{eq:QRbound}
\eea
\end{conjecture}
It is immediately clear that (\ref{eq:QRbound}) implies (\ref{eq:the2}),
as can be seen by setting
$Q:=[A_j^* A_k]$, $R:=[B_j^* B_k]$,
$q_k=a_k^2$ and $r_k=b_k^2$.
When $N=2$, the statement of Conjecture \ref{conj2} coincides with Conjecture \ref{conj3} and is therefore equivalent with it.
For larger $N$ we don't know whether equivalence holds between Conjectures \ref{conj1} and \ref{conj2}.

We can further reformulate Conjecture \ref{conj2} once we realise that it is true for matrices $Q$ and $R$
consisting of scalar blocks, which is proven in Section \ref{sec:rank1}.
Namely, Conjecture \ref{conj2} is thus equivalent to the statement that the maximum
in the LHS of (\ref{eq:QRbound}) is obtained when the blocks of $Q$ and $R$ are essentially scalars
(i.e.\ all blocks are proportional to the same matrix).

\bigskip

In an entirely similar way, we can treat the case $1\le p\le 2$, which corresponds to $1/2\le q\le1$.
Actually, the equivalences hold in the more general case $0\le q\le 1$, although the statements themselves
are no longer valid for $q\le 1/2$.
In the above, we just have to literally replace `maximisation' by `minimisation', `convex' by `concave', and
reverse the inequalities.
\section{Proofs in Special Cases}
\subsection{The norm compression has rank 1.}
Let $A_k$ and $B_k$ be such that $a_k^2 = \alpha^2 p_k$ and $b_k^2 = \beta^2 p_k$, with $p_k\ge0$
and $\sum_k p_k=1$. This is indeed what is meant with $T$ having a rank 1 norm compression.
The choice of $p_k$ is such that they form a probability distribution.
By setting $A_k=\sqrt{p_k}\,X_k$ and $B_k=\sqrt{p_k}\,Y_k$,
with $||X_k X_k^*||_q=\alpha^2$ and $||Y_k Y_k^*||_q =\beta^2$,
(\ref{eq:the2}) reduces to
$$
\schatten{q}{
\sum_k p_k \twomat{X_k X_k^*}{\quad X_k Y_k^*}{Y_k X_k^*}{\quad Y_k Y_k^*}
} \le
\schatten{q}{\twomat{\alpha^2}{\quad\alpha\beta}{\alpha\beta}{\quad\beta^2}},
$$
for $q\ge 1$, and reversed for $0\le q\le1$.
Now note that the RHS is independent of $p_k$.
By convexity of the Schatten norms for $q\ge1$, and concavity of the Schatten quasi-norms for $0\le q\le 1$,
it is therefore enough to prove this inequality for the
extremal points $(p_k)_k = (0,\ldots,0,1,0,\ldots,0)$, which amounts to the case that $N=1$.
That is, we need to prove
$$
\schatten{q}{
\twomat{X_1 X_1^*}{\quad X_1 Y_1^*}{Y_1 X_1^*}{\quad Y_1 Y_1^*}
} \le \schatten{q}{
\twomat{\alpha^2}{\quad \alpha\beta}{\alpha\beta}{\quad \beta^2}
}.
$$
As the left-hand side is equal to $\schatten{q}{X_1^*X_1+Y_1^*Y_1}$, and the right-hand side
is equal to $\alpha^2+\beta^2=||X_1^* X_1||_q+||Y_1^* Y_1||_q$, this simply follows from the triangle inequality;
that is, again from convexity of the Schatten norms for $q\ge 1$.
The reversed inequality for $0\le q\le1$ follows similarly from concavity of the Schatten $q$-quasi-norms on PSD matrices.
\qed
\subsection{All $A_k$ and $B_k$ have rank 1.\label{sec:rank1}}
In this case, $A_k$ and $B_k$ can be written as
$$
A_k = \alpha_k u_k^*,\quad B_k = \beta_k v_k^*,
$$
with $\alpha_k$, $\beta_k$, $u_k$ and $v_k$ vectors, such that
$||\alpha_k||=a_k$, $||\beta_k||=b_k$, $||u_k||=1$ and $||v_k||=1$.
Then $||A_k||_p=a_k$ and $||B_k||_p=b_k$ for all $p>0$.

The left-hand side of (\ref{eq:the2}) then can be written as
$$
\LHS(\ref{eq:the2})=\schatten{q}{\twomat{\sum_k \alpha_k \alpha_k^*}{\sum_k (u_k^* v_k)\alpha_k\beta_k^*}
{\sum_k (v_k^* u_k)\beta_k\alpha_k^*}{\sum_k \beta_k\beta_k^*}}.
$$
Introducing $K=\oplus_k u_k^* v_k$, and denoting by $\cA$ and $\cB$ the matrices whose columns are $\alpha_k$ and $\beta_k$, respectively,
this is equal to
$$
\LHS(\ref{eq:the2})=\schatten{q}{\twomat{\cA\cA^*}{\cA K \cB^*}{\cB K^* \cA^*}{\cB\cB^*}}.
$$
The $q$-th power of the LHS equals the trace of the $q$-th power of the given block matrix, which is PSD.
Therefore, for $q\ge1$, the $q$-th power of the LHS
is convex in $K$ and $K^*$, while it is concave in it for $0\le q\le1$.
Now note that $K$ is a diagonal contraction. Since the diagonal contractions
form a convex set with extremal points the diagonal unitaries, the LHS is extremal for $K$ equal to some diagonal unitary $U$.
Therefore, for $q\ge1$, introducing the matrices $Q=\cA^*\cA$ and $R=U\cB^*\cB U^*$,
$$
\LHS(\ref{eq:the2}) \le \schatten{q}{\twomat{\cA\cA^*}{\cA U \cB^*}{\cB U^* \cA^*}{\cB\cB^*}}
= \schatten{q}{\cA^*\cA + U\cB^*\cB U^*}
= \schatten{q}{Q+R},
$$
whereas, for $0\le q\le 1$, the direction of the inequality is reversed.
Note that the diagonal entries of $Q$ and $R$
are given by $Q_{kk} = (\cA^*\cA)_{kk}$ and, since $U$ is a \textit{diagonal} unitary, $R_{kk} = (\cB^*\cB)_{kk}$.

The scalars $a_k$ and $b_k$ appearing in the right-hand side of (\ref{eq:the2}) are thus given by
$a_k = ||\alpha_k|| = \sqrt{(\cA^*\cA)_{kk}} = \sqrt{Q_{kk}}$, and
$b_k = ||\beta_k|| = \sqrt{(\cB^*\cB)_{kk}} = \sqrt{R_{kk}}$.
Hence, the right-hand side of (\ref{eq:the2}) depends only on the diagonal elements of $Q$
and $R$. Defining $q_k:=Q_{kk}$ and $r_k:=R_{kk}$, it is given by
$$
\RHS(\ref{eq:the2}) = \schatten{q}{
\twomat{\sum_k q_k}{\sum_k\sqrt{q_k r_k}}{\sum_k\sqrt{q_k r_k}}{\sum_k r_k}
} =
\schatten{q}{
\twomat{\trace Q}{y}{y}{\trace R}
},
$$
with $y:=\sum_k\sqrt{q_k r_k}$.
In other words, the special case under consideration follows from the special case of Conjecture \ref{conj2} where
the blocks of $Q$ and $R$ are scalar, and which I will now prove.

In fact, it holds for any unitarily invariant (quasi)norm, and not just for the Schatten norms.
More precisely,
\be\label{eq:the4}
|||Q+R||| \le \left|\left|\left|\,\twomat{\trace Q}{y}{y}{\trace R}\oplus 0_{d-2}\,\right|\right|\right|.
\ee
Furthermore, for $0\le q\le 1$, the reversed inequality holds for the Schatten quasi-norms.
To prove this, I only have to show (\ref{eq:the4}) for the Ky Fan norms $||.||_{(k)}$, $k=1,2,\ldots,d$,
because I can then invoke Ky Fan's Dominance Theorem.
This Theorem is formulated for norms, but can easily be generalised to concave quasi-norms.
If a Hermitian matrix $A$ dominates another Hermitian matrix $B$ in the sense that its eigenvalues majorise those of the other,
$\lambda(A)\succ\lambda(B)$, then $A$ will be dominated by $B$ in any unitarily invariant concave quasi-norm.
In particular, it will be dominated in any Schatten $q$-quasi-norm with $0\le q\le 1$.
This follows directly from the statement that $x\prec y$ implies $\sum_{j=1}^d \phi(x_j)\ge\sum_{j=1}^d\phi(y_j)$
for concave functions $\phi$ (see, e.g.\ Theorem II.3.1 in \cite{Bhatia} mutatis mutandis for concave $\phi$).

Thus we only have to prove (\ref{eq:the4}) for the Ky Fan norms. We need only consider two cases.
The case $k>1$ is trivial, as the right-hand side is then nothing but the trace of the matrix.
This trace is $\trace(Q+R)$, which is obviously an upper bound on $||Q+R||_{(k)}$.

The remaining case $k=1$, i.e.\ the operator norm, or $q=\infty$,
is only slightly more complicated. I will prove (\ref{eq:the4}) for the operator norm
by showing that the right-hand side is the maximum of $||Q+R||_\infty$ over all $Q,R\ge0$ with prescribed diagonals
$Q_{kk}=q_k$, $R_{kk}=r_k$. In fact, the following Lemma is even a little more general:
\begin{lemma}
Let $q_k$, $r_k$, $s_k$, \ldots, with $k=1,\ldots,d$, be given non-negative numbers.
The maximal value of $||Q+R+S+\ldots||_\infty$ over all $d\times d$ PSD matrices $Q$, $R$, $S$,\ldots
with prescribed diagonal elements $q_k$, $r_k$, $s_k$, \ldots, respectively,
is obtained when $Q$, $R$, $S$,\ldots are rank 1 matrices with non-negative elements.
That is, $Q_{ij}=\sqrt{q_i q_j}$, etc.
\end{lemma}
\textit{Proof.}
Since $Q$, $R$, $S$,\ldots are Hermitian, the norm $||Q+R+S+\ldots||_\infty$ is given by
$\max_\psi \psi^*(Q+R+S+\ldots)\psi = \max_\psi \psi^*Q\psi+\psi^*R\psi+\ldots$, where $\psi$ is any normalised vector.
Now, for any vector $\psi$ we have the inequality
\beas
\psi^* Q \psi &=& \sum_k |\psi_k|^2 q_k + 2\sum_{j<k}\re( \overline{\psi}_j \psi_k Q_{jk}) \\
&\le& \sum_k |\psi_k|^2 q_k + 2\sum_{j<k} |\psi_j| |\psi_k| \sqrt{q_j q_k},
\eeas
where $|Q_{jk}|\le \sqrt{q_j q_k}$ is required for positivity of $Q$.
Letting $\psi_j=\exp(i \phi_j)|\psi_j|$, equality is obtained when $Q_{jk}=\exp(i(\phi_j-\phi_k))\sqrt{q_j q_k}$.
Since the maximal value of $\psi^*Q\psi$ for this $Q$ does not depend on the arguments $\phi_j$,
we could as well restrict $\psi$ to positive real vectors.
As the above is true for any value of $\psi$, $||Q+R+S+\ldots||_\infty$ is maximal for $Q_{jk}=\sqrt{q_j q_k}$,
$R_{jk}=\sqrt{r_j r_k}$, \ldots
\qed

With this Lemma, we are done.\qed

\textit{Remark:} The Lemma cannot be generalised to other norms, because the extremal points of the set of positive matrices
with prescribed diagonal entries do not necessarily have rank 1, even in the simplest case that these diagonal entries are all 1
(the case of so-called correlation matrices \cite{ckli}).
\subsection{All $A_k$ are proportional to each other, and so are all $B_k$}
Let us now consider the case where the $A_k$ satisfy $A_k = \alpha_k X$, for some scalars $\alpha_k$
and some matrix $X$, and similarly, $B_k=\beta_k Y$.
For this special case we need a result by King from \cite{kingEB}:
\begin{proposition}[King]
For $A_k, B_k\ge 0$, and any $q\ge 1$,
\be
\schatten{q}{\sum_k A_k\otimes B_k} \le \schatten{q}{\sum_k A_k} \, \max_j \schatten{q}{B_j},
\ee
while for $0\le q\le 1$,
\be
\schatten{q}{\sum_k A_k\otimes B_k} \ge \schatten{q}{\sum_k A_k} \, \min_j \schatten{q}{B_j}.
\ee
\end{proposition}
\textit{Proof.}
Since the $A_k$ are positive, the following notations are meaningful:
\beas
F   &=& (\sqrt{A_1} \otimes \identity \ldots \sqrt{A_K} \otimes \identity) \\
G   &=& (\sqrt{A_1} \ldots \sqrt{A_K}) \\
H   &=& \bigoplus_k \identity\otimes B_k.
\eeas
Let us denote by $X_{(kk)}$ the $k$-th diagonal block of a matrix in the
same partitioning as $H$.
For example, $H_{(kk)}=\identity\otimes B_k$.

Using these notations, $\sum_k A_k\otimes B_k$ can be written as $FHF^*$.
The Araki-Lieb-Thirring inequality (\cite{Bhatia}, Theorem IX.2.10) can be easily brought in a form
that states that, for $H\ge0$ and general $F$,
$\trace[(FHF^*)^q]\le \trace[(F^* F)^q \, H^q]$ holds for $q\ge1$, and the reversed inequality
for $0\le q\le1$. Thus, for $q\ge1$,
\beas
\schatten{q}{\sum_k A_k\otimes B_k}^q &=& \trace[(FHF^*)^q] \\
&\le& \trace[(F^* F)^q \, H^q] \\
&=& \sum_k \trace[[(F^* F)^q]_{(kk)} \, (\id \otimes B_k^q)] \\
&=& \sum_k \trace[[(G^* G)^q]_{(kk)}] \, \trace[B_k^q] \\
&\le& \max_j \trace[B_j^q] \, \sum_k \trace[[(G^* G)^q]_{(kk)}] \\
&=& \max_j \trace[B_j^q] \, \trace[(G^* G)^q].
\eeas
while for $0\le q\le 1$, the direction of the inequalities reverse, and the `max' has to be replaced
by a `min'.
Then noting
$$
\trace[(G^* G)^q] =  \trace[(G G^*)^q] =  \trace[(\sum_k A_k)^q],
$$
and taking $q$-th roots yields the Proposition.
\qed

This proposition can be reformulated as:
\begin{corollary}
For $A_k,B_k\ge 0$, and any $q\ge 1$,
\be
\schatten{q}{\sum_k A_k\otimes B_k} \le \schatten{q}{\sum_k A_k \schatten{q}{B_k}},
\ee
while for $0\le q\le 1$ the inequality reverses.
\end{corollary}
\textit{Proof.}
Define $A'_k = ||B_k||_q A_k$ and $B'_k = B_k/||B_k||_q$.
Then $||B'_k||_q=1$ and, for $q\ge1$,
\beas
\schatten{q}{\sum_k A_k\otimes B_k} &=& \schatten{q}{\sum_k A'_k\otimes B'_k}
\le \max_j \schatten{q}{B'_j} \, \schatten{q}{\sum_k A'_k} \\
&=& \schatten{q}{\sum_k ||B_k||_q A_k}.
\eeas
The proof is completely similar for $0\le q\le 1$.
\qed

What we need now is:
\begin{corollary}\label{co:king2}
For general matrices $X_k$, $k=1,\ldots,d$, with $||X_k X_k^*||_q=1$ and
$Z=[z_{ij}]_{i,j=1}^d\ge0$,
the block matrix $Q:=[z_{ij} X_i X_j^*]$
satisfies $||Q||_q \le ||Z||_q$ for $q\ge1$ and $||Q||_q \ge ||Z||_q$ for $0\le q\le 1$.
\end{corollary}
\textit{Proof.}
We can write
$$
Q = (\oplus_i X_i) \,\,(Z\otimes\id)\,\,(\oplus_i X_i^*),
$$
thus $Q$ is unitarily equivalent with
$$
(Z^{1/2}\otimes\id)\,\,(\oplus_i X_i X_i^*) \,\,(Z^{1/2}\otimes\id)
$$
and has the same $q$-norm.
Denoting the columns of $Z^{1/2}$ by $y_i$, the latter matrix can be written as
$$
\sum_i y_i y_i^*\otimes X_i X_i^*.
$$
Thus, $||Q||_q = ||\sum_i y_i y_i^*\otimes X_i X_i^*||_q$, which
by the previous Corollary is upper bounded by $||\sum_i y_i y_i^* ||X_i X_i^*||_q\,||_q$, for $q\ge1$.
By the assumption $||X_i X_i^*||_q=1$, this upper bound reduces to $||\sum_i y_i y_i^*||_q = ||Z||_q$.

For $0\le q\le 1$, the proof is completely similar; replace `upper bound' by `lower bound'.
\qed

Now put $A_k=a_k X$, $B_k=b_k Y$, with $X$ and $Y$ such that $||XX^*||_q=||YY^*||_q=1$.
With this choice the left-hand side of (\ref{eq:the2}) becomes
$$
\schatten{q}{
\twomat{(\sum_k a_k^2)\, XX^*}{\quad(\sum_k a_k b_k)\, XY^*}{(\sum_k a_k b_k)\, YX^*}{\quad(\sum_k b_k^2)\, YY^*}
},
$$
and the validity of (\ref{eq:the2}) in this case follows directly from Corollary \ref{co:king2},
with
$$
Z = \twomat{\sum_k a_k^2}{\quad\sum_k a_k b_k}{\sum_k a_k b_k}{\quad\sum_k b_k^2} \ge0.
$$
\subsection{The case $q\ge2$}
For this case, which corresponds to $p\ge 4$, I start with a simple Lemma about $2\times 2$ PSD matrices; the positivity
is essential.
\begin{lemma}\label{lem:mono}
Let $A_i = \twomat{a_i}{c_i}{c_i}{b_i}$ be $2\times2$ PSD matrices with non-negative elements.
For any UI norm $|||.|||$, if $a_1\le a_2$, $b_1\le b_2$ and $c_1\le c_2$, then $|||A_1||| \le |||A_2|||$.
\end{lemma}
\textit{Proof.}
By Ky Fan's dominance theorem \cite{HJII,Bhatia}, one only needs to show the statement
for the trace norm and operator norm. Since $A_1$ and $A_2$
are PSD, their trace norm equals their trace, and obviously $\trace A_1\le \trace A_2$.
For the operator norm, note that
$$
||A_1||_\infty = \max_\psi \psi^* A_1 \psi = a_1|\psi_1|^2 + b_1|\psi_2|^2+2c_1\re\overline{\psi}_1\psi_2.
$$
Since $c_1\ge0$, in order for the maximum to be achieved, $\re\overline{\psi}_1\psi_2$ must be positive.
In that case one sees that $\psi^* A_1\psi \le \psi^* A_2\psi$, so the same must hold for the respective
maximums.
\qed

From King's norm compression inequality for PSD $2\times2$ block matrices mentioned in the Introduction
(i.e.\ (\ref{eq:the}) with $N=2$ and $T\ge0$), it follows that, for $q\ge2$,
$$
\schatten{q}{
\twomat{\sum_k A_k A_k^*}{\quad\sum_k A_k B_k^*}{\sum_k B_k A_k^*}{\quad\sum_k B_k B_k^*}
}
\le
\schatten{q}{
\twomat{||\sum_k A_k A_k^*||_q}{\quad||\sum_k A_k B_k^*||_q}{||\sum_k B_k A_k^*||_q}{\quad||\sum_k B_k B_k^*||_q}
}.
$$
By the triangle inequality, and the Cauchy-Schwarz inequality for matrices (\cite{Bhatia} IX.30),
the matrix elements of the right-hand side are upper bounded as follows:
\beas
||\sum_k A_k A_k^*||_q &\le& \sum_k ||A_k A_k^*||_q = \sum_k a_k^2 \\
||\sum_k B_k B_k^*||_q &\le& \sum_k b_k^2 \\
||\sum_k A_k B_k^*||_q &\le& \sum_k ||A_k B_k^*||_q \le \sum_k ||A_k A_k^*||_q^{1/2} ||B_k B_k^*||_q^{1/2} = \sum_k a_k b_k.
\eeas
Application of Lemma \ref{lem:mono} then immediately proves (\ref{eq:the2}) for $q\ge 2$.\qed

\textit{Remarks.}
\begin{enumerate}
\item For $1\le q\le2$ one might be tempted to apply the NCI proved in \cite{ka} instead of King's, but, unfortunately, it
gives a bound which can be higher than (\ref{eq:the2}).
\item Note also that, since we don't have the required duality result, this case does not allow
us to treat the $1\le p\le 4/3$ case. The only possibility to do so would be to
use a variant of King's inequality for $q\le 1$, which we don't have.
Preliminary numerical experiments indicate that for $0\le q\le 1$, and $T$ a $2\times2$ PSD block matrix
the following holds:
$$
\schatten{q}{T}\le \schatten{q}{\cC_q(T)},
$$
just as in the $q\ge2$ case, which is in the wrong direction.
\item Similarly, the NCI proved in \cite{ka} also seems to work for $0\le q\le 1$, with the same ordering
as for $q\ge2$.
\end{enumerate}
\section{Diagonal blocks, the $0\le p\le 1$ case, and a counterexample for $3\times 3$ block matrices}
In this final Section, we study the special case where all blocks are diagonal.
For $N=2$, this case has been considered and proven in \cite{king_nath}.
The proof proceeds via reducing the problem to the following convexity result:
\begin{lemma}\label{lem:king_nath}
Define the (homogeneous) real-valued function
$$
B\mapsto g_{p}(B) = \schatten{p}{B^{\circ 1/p}}^{p}
$$
on $M_{2,2}(\R^+)$,
where $B^{\circ 1/p}$ denotes the entry-wise (Hadamard) power of $B$.
Then $g_{p}(B)$ is convex for $1\le p\le 2$, while it is concave for $2\le p$.
\end{lemma}
Numerical experiments quickly reveal that for $0\le p\le 1$, $g_{p}(B)$ is neither convex nor concave.

With this Lemma at hand, it is straightforward to prove (\ref{eq:the})
for diagonal blocks. The question thus arises whether the Lemma remains true for matrices in $M_{2,N}(\R^+)$,
or even for general non-negative matrices.

\textit{Proof that (\ref{eq:the}) for diagonal blocks is implied by Lemma \ref{lem:king_nath}.}
We present the proof in the most general setting of matrices of arbitrary size, although we will show
below that the Lemma does not hold for matrices of size $3\times3$ (and larger).

Consider the case $p\ge2$.

Let the block matrix under consideration be denoted $T = [T_{(ij)}]$ and
let $a_{ij}^k$ denote the $k$-th diagonal entry of $T_{(ij)}$.
With these diagonal entries we can construct
$d$ matrices $[a_{ij}^k]_{ij}$.
By a simple unitary conjugation, $T$ can thus be brought into the direct sum form
$\bigoplus_{k=1}^d [a_{ij}^k]_{ij}$,
so that
$$
||T||_p^p = \sum_k \schatten{p}{[a_{ij}^k]_{ij}}^p.
$$
Every matrix in this sum can be trivially norm compressed by applying the absolute value
to every matrix entry.
For this norm compression, the NCI (\ref{eq:the}) applies, because this is a special instance of the case where
all blocks are rank 1.
For $p\ge 2$, this gives
$$
||T||_p^p \le \sum_j \schatten{p}{[|a_{ij}^k|]_{ij}}^p.
$$
Since all matrix entries are now non-negative, Lemma \ref{lem:king_nath}
can be applied, giving
$$
||T||_p^p \le \schatten{p}{
\left[(\sum_k |a_{ij}^k|^p)^{1/p}\right]_{ij}
}^p
=\schatten{p}{[ ||T_{(ij)}||_p ]}^p.
$$
One observes that the right-hand side is just the right-hand side of (\ref{eq:the}).

For $1\le p\le 2$, we proceed in a completely similar way.
\qed

Upon closer inspection, we see that validity of Conjecture \ref{conj1}
in the case of \textit{positive} diagonal blocks is equivalent to validity of Lemma \ref{lem:king_nath}.
Hence, counterexamples to the conjectured validity of Lemma in the $2\times N$ case yield counterexamples to Conjecture \ref{conj1}.

For $0\le p\le 1$, we already noted that $g_{p}(B)$ is neither convex nor concave,
whence (\ref{eq:thedual}) is not true for $0\le p\le 1$. This is exhibited
by $T$ consisting of positive $2\times 2$ diagonal blocks, constructed from two matrices in $M_2(\R^+)$ that
violate convexity of $g_{p}(B)$ (i.e.\ from $A$ and $B$ such that $g_p(A+B)\not\le g(A)+g(B)$).

Proving the Lemma for $2\times N$ matrices turns out to be surprisingly hard. The proof of the $2\times 2$ case in \cite{king_nath}
is already very involved and no simple method seems to be forthcoming yet.
In the remainder of this Section we introduce a certain direction along which a proof may eventually
be found; using this method we do find a counterexample for Conjecture \ref{conj1} in the $3\times3$ case.

\bigskip

Let us focus attention to the function $g_p(B)$ defined above and to the case $1\le p\le 2$.
Following \cite{king_nath} we express the convexity using a differential.
By homogeneity of $g_p$, its convexity is equal to its subadditivity; i.e.\ we have to show that for $B,\Delta\in M(\R^+)$,
$g_p(B+\Delta) \le g_p(B)+g_p(\Delta)$. Replacing $\Delta$ by $t\Delta$, $t>0$, this yields
$$
g_p(B+t\Delta) \le g_p(B)+t g_p(\Delta).
$$
For infinitesimal $t$ we get the requirement
$$
\ddt g_p(B+t\Delta) \le g_p(\Delta).
$$
This is actually an equivalent statement to the former one, as can be seen by integrating over $t$.

The derivative can be calculated explicitly.
For that purpose we need to know the Fr\'echet derivative of a Schatten norm $\ddt \trace|B+t\Delta|^p$.
For positive $B$ and $\Delta$, the absolute value is not needed and we simply get
$$
\ddt \trace(B+t\Delta)^p = p\trace(B^{p-1}\Delta).
$$
In the general case one can show that the $(p-1)$-th power of $B$ has to be replaced by the quantity
$(BB^*)^{p/2-1}B$. In terms of the polar decomposition $B=U|B|$, this quantity is equal to $|B^*|^{p-1}U$,
which is why I propose to call this quantity the $(p-1)$-th \textit{polar power} of $B$.
I will use the shorthand $B^{[r]}$ for the $r$-th polar power:
$$
B^{[r]} := (BB^*)^{(r-1)/2}B.
$$
Some obvious statements one can make are that
the polar power coincides with the ordinary power on PSD matrices,
and that UI norms don't distinguish between polar powers and ordinary powers: $|||B^{[r]}||| = |||\,\,|B|^r\,\, |||$.

With the derivative of the Schatten norm at hand, the Fr\'echet derivative of $g_p(B)$ can be calculated,
and the result of a somewhat lengthy calculation is
$$
\ddt g_p(B+t\Delta) = \trace \Delta^T((B^{\circ 1/p})^{[p-1]} \circ B^{\circ (1-p)/p}).
$$
Recall that $\circ$ denotes the Hadamard product, and $B^{\circ r}$ denotes the $r$-th entry-wise power.
With the substitutions $A:=\Delta^{\circ 1/p}$ and $C:=B^{\circ 1/p}$, the statement we have to prove
becomes
\be\label{eq:AC}
\trace (A^{\circ p})^T(C^{[p-1]}\circ C^{\circ 1-p}) \le \trace|A|^p,
\ee
for $1\le p\le 2$.
Note that equality holds when $A$ and $C$ are proportional.

\bigskip

I have not been able to prove inequality (\ref{eq:AC}) for $A,C\in M_{2,N}(\R^+)$. However, by studying what happens
when $C$ is on the boundary of $M(\R^+)$ (that is, when some entries of $C$ are 0) I have been able to
find a counterexample for $M_3(\R^+)$, which, as mentioned above, automatically yields a counterexample to
Conjecture \ref{conj1} for $3\times 3$ block matrices.

Consider $p=1.5$, say, and the matrix
$$
C=\left(
\begin{array}{ccc}
\epsilon &1& \epsilon \\
1&\epsilon&1 \\
\epsilon&1&1
\end{array}
\right),
$$
where $\epsilon$ will tend to 0.
One can easily calculate the polar power of $C$ numerically. We will only need entry $1,1$ to see the violation
of the inequality here. The result is
$$
\lim_{\epsilon\to 0}(C^{[0.5]})_{1,1} = 0.11669...,
$$
where the only thing that matters here is that this entry is strictly positive.
Indeed, the $1,1$ entry of $C^{\circ 1-p}$ is given by $\epsilon^{-0.5}$, which tends to $+\infty$ as $\epsilon$
tends to 0, and, therefore, the $1,1$ entry of $(C^{[p-1]}\circ C^{\circ 1-p})$ tends to $+\infty$ as well.
Hence, for any $A$ with non-zero $1,1$ entry, the inequality (\ref{eq:AC}) is violated to the maximally possible extent.

\bigskip

In the light of this counterexample, one may get concerned about validity of (\ref{eq:AC}) in the $2\times N$ case
as well. I end this Section by showing that the phenomenon just mentioned,
of certain entries of $(C^{[p-1]}\circ C^{\circ 1-p})$ tending to $+\infty$, will not occur in $2\times N$ matrices
(but leaving open the general question of validity of (\ref{eq:AC})).
Indeed, we have the following Lemma:
\begin{lemma}
Let $C$ be a bounded matrix in $M_{2,N}(\R^+)$, $1\le p\le 2$.
If entry $C_{ij}=0$, then the corresponding entry $(C^{[p-1]})_{ij}$ is non-positive.
\end{lemma}
It is, of course, this non-positivity that rescues inequality (\ref{eq:AC}) here.
If entry $C_{ij}$ is zero, entry $(C^{\circ 1-p})_{ij}$ tends to $+\infty$, so that
$(C^{[p-1]}\circ C^{\circ 1-p})_{ij}$ is either 0 or $-\infty$. In the former case, it does not contribute to the LHS
of (\ref{eq:AC}), while in the latter it is a dominating contribution (at least, for $A$ with non-zero entries)
and causes (\ref{eq:AC}) to be satisfied irrespective of the value of any other entry.

\textit{Proof.}
By definition, $C^{[p-1]} = |C^T|^{p-2}C$. Here, $p-2<0$.
Now, $|C^T|$ is a PSD $2\times 2$ matrix with non-negative entries.

This is also true for any positive power of $|C^T|$.
Indeed, a useful parametrisation of PSD $2\times 2$ matrices with non-negative entries is
$A=\lambda P+\mu (\id-P)$, where $\lambda\ge\mu$ are
the eigenvalues of $A$, and the projector $P$ is given as
$P=\twomat{t}{s}{s}{1-t}$, with $0\le t\le 1$ and $s=\sqrt{t(1-t)}\in[0,1/2]$.
Obviously, for any $p$, $A^p=\lambda^p P+\mu^p (\id-P)$.
Written out,
$$
A^p = \twomat{t\lambda^p + (1-t)\mu^p}{s(\lambda^p-\mu^p)}{s(\lambda^p-\mu^p)}{(1-t)\lambda^p + t\mu^p},
$$
which is clearly entry-wise non-negative for all $p\ge0$.
For negative $p$, on the other hand, we see that $(A^p)_{12}$ reverses sign (or stays 0).

Since $p-2$ is indeed negative, we thus find that $(|C^T|^{p-2})_{12}$ is non-positive.
Consider then the zero entry $C_{1j}=0$ (the reasoning is identical for row 2). Then
$$
(C^{[p-1]})_{1j} = (|C^T|^{p-2})_{11} C_{1j} + (|C^T|^{p-2})_{12} C_{2j} = (|C^T|^{p-2})_{12} C_{2j},
$$
which is non-positive.
\qed

\begin{ack}
This work was supported by The Leverhulme Trust (grant F/07 058/U),
and is part of the QIP-IRC (www.qipirc.org) supported by EPSRC (GR/S82176/0).
The author acknowledges support from the Institute for Mathematical Sciences, Imperial College London.
Sincere thanks to Chris King for finding yet another duality mistake in a previous version of the manuscript.
\end{ack}


\begin{thebibliography}{99}
\bibitem{Bhatia} R.~Bhatia, \textit{Matrix Analysis}, Springer, Heidelberg (1997).
\bibitem{carlenlieb} E.A.~Carlen and E.H.~Lieb,
Advances in Math.\ Sciences, AMS Transl.\ (2) \textbf{189}, 59--62 (1999). See also arxiv.org E-print math.OA/0701352.
\bibitem{HJII} R.A.~Horn and C.R.~Johnson, \textit{Topics in Matrix Analysis},
Cambridge University Press, Cambridge (1991).
\bibitem{lieb} K.~Ball, E.~Carlen and E.H.~Lieb, ``Sharp uniform convexity and smoothness
inequalities for trace norms,'' Invent.\ Math.\ \textbf{115}, 463--482 (1994).
\bibitem{ka} K.M.R.~Audenaert, ``A norm compression inequality for block partitioned positive definite matrices,''
Lin.\ Alg.\ Appl.\ \textbf{413}, 155--176 (2006).
\bibitem{kingEB} C.~King, ``Maximal $p$-norms of entanglement-breaking channels,'' Quantum Information and Computation \textbf{3}, 186--190 (2003).
\bibitem{king} C.~King, ``Inequalities for trace norms of $2\times 2$ block matrices,'' Commun.\ Math.\ Phys.\ \textbf{242}, 531--545 (2003).
\bibitem{king_nath} C.~King and M.~Nathanson, ``New trace norm inequalities for $2\times 2$ blocks of diagonal matrices,''
Lin.\ Alg.\ Appl.\ \textbf{389}, 77--93 (2004).
\bibitem{hanner} O.~Hanner, ``On the uniform convexity of $L^p$ and $l^p$,'' Ark.\ Math.\ \textbf{3}, 239--244 (1958).
\bibitem{TJ} N.~Tomczak-Jaegermann, ``The moduli of smoothness and convexity and Rademacher averages of trace classes $S_p$,''
Studia Math.\ \textbf{50}, 163--182 (1974).
\bibitem{ckli} C.K.~Li and B-S. Tam, ``A note on extreme correlation matrices,'' SIAM J.\ Matrix Anal.\ Appl.\ \textbf{15}, 903--908 (1994).
\end{thebibliography}
\end{document}